\documentclass[a4paper]{article}
\usepackage[T1]{fontenc}
\usepackage{latexsym}
\usepackage{amssymb,amsmath}

\newtheorem{Theorem} {Theorem} [section]
\newtheorem{Proposition} [Theorem] {Proposition}
\newtheorem{Corollary}[Theorem]{Corollary}
\newcommand{\Proof}{ \noindent{\bf Proof.}\quad }
\newcommand{\qed}{\hfill$\Box$\medskip}
\newcommand{\Ff}{{\mathbb F}}
\newcommand{\Zz}{{\mathbb Z}}
\newcommand{\Cc}{{\mathbb C}}
\newcommand{\adj}{{~\sim~}}

\newcommand{\tin}[6]{[\begin{smallmatrix}
#1 & #2 & #3 \\
#4 & #5 & #6 \\[1pt]
\end{smallmatrix}]}
\newcommand{\tinA}[2]{[\begin{smallmatrix}
#1 \\
#2
\end{smallmatrix}]}
\newcommand{\tinB}[4]{[\begin{smallmatrix}
#1 & #2 \\
#3 & #4
\end{smallmatrix}]}
\newcommand{\tinC}[6]{[\begin{smallmatrix}
#1 & #2 & #3 \\
#4 & #5 & #6
\end{smallmatrix}]}
\newcommand{\tinD}[8]{[\begin{smallmatrix}
#1 & #2 & #3 & #4 \\
#5 & #6 & #7 & #8
\end{smallmatrix}]}

\title{Triple intersection numbers for the Paley graphs}
\author{Andries E.~Brouwer \\ {\tt aeb@cwi.nl}
\and William J.~Martin\footnote{
Worcester Polytechnic Institute,
Dept. of Mathematical Sciences, Worcester, MA USA} \\
{\tt martin@wpi.edu}}

\date{2021-09-05}

\begin{document}
\maketitle

\vspace{-0.1cm}
\begin{abstract}
We give a tight bound for the triple intersection numbers
of Paley~graphs. In particular, we show that any three vertices
have a common neighbor in Paley~graphs of order larger than 25.
\end{abstract}

\bigskip
Let $q = 4t+1$ be a prime power, and let $\Gamma$ be Paley($q$),
the Paley graph on $q$ vertices, with as vertex set the finite field
$\Ff_q$ of size $q$, where two vertices are adjacent when their difference
belongs to $\Ff_q^{*2}$, the set of  nonzero squares in $\Ff_q$.
This graph is connected with diameter 2, and self-complementary. 

In \cite{MartinWashock}, the authors needed the fact that any function
$\psi: \Ff_q^{*2} \cup \{0\} \rightarrow \Cc^*$ satisfying 
{\it (i)} $\psi(0)=1$ and
{\it (ii)} $\psi(a)\psi(b)=\psi(c)\psi(d)$ whenever $a+b=c+d$ 
must be the restriction of some additive character of $\Ff_q$ if $q > 5$.
The present note provides a short proof of that fact.
 
Following the notation of \cite{CoolsaetJurisic08} \S 3, define
\emph{generalized intersection numbers}
$\tinD{a_1}{a_2}{\cdots}{a_\ell}{i_1}{i_2}{\cdots}{i_\ell}$
for $a_1,\ldots, a_\ell \in \Ff_q$ and $i_1,\ldots,i_\ell \in \{0,1,2\}$
by $\tinD{a_1}{a_2}{\cdots}{a_\ell}{i_1}{i_2}{\cdots}{i_\ell} :=
| \Gamma_{i_1}(a_1) \cap \cdots \cap \Gamma_{i_\ell}(a_\ell) |$, where 
$\Gamma_i(a)$ denotes the set of vertices at distance $i$ from $a$.
Note that $\sum_{i_\ell} \tinC{a_1}{\cdots}{a_\ell}{i_1}{\cdots}{i_\ell}
=  \tinC{a_1}{\cdots}{a_{\smash{\ell-1}}}{i_1}{\cdots}{i_{\smash{\ell-1}}}$
and $\tinA{a}{i} = \frac{q-1}{2}$
and $\tinB{a}{b}{i}{j} = \frac{q-1}{4} - \delta_{hi} \delta_{hj} \delta_{ij}$
for distinct $a,b$ and $h,i,j=1,2$ 
where $h$ is the distance from $a$ to $b$.
It follows that all $\tin{a}{b}{c}{h}{i}{j}$ are known if one knows
$\tin{a}{b}{c}{1}{1}{1}$.

\begin{Proposition}\label{bds}
$\left|  \tin{a}{b}{c}{1}{1}{1} - \frac{q-9}{8} \right| \le
\frac14 \sqrt{q} + \frac34$ for any three distinct $a,b,c$.
\end{Proposition}

\Proof
Let $\chi$ be the quadratic character.
If $a,b,c$ are distinct, then
$$
\smash{\sum_x} (1+\chi(x-a))(1+\chi(x-b))(1+\chi(x-c)) =
8 \,\tin{a}{b}{c}{1}{1}{1}  + 4 R
$$
where $R = \tin{a}{b}{c}{0}{1}{1} + \tin{a}{b}{c}{1}{0}{1} +
\tin{a}{b}{c}{1}{1}{0}$, so that $R \in \{0,1,3\}$.
Let $S = \sum_x \chi((x-a)(x-b)(x-c))$.
Since $\sum_x 1 = q$ and $\sum_x \chi(x) = 0$
and $\sum_x \chi(x(x-a)) = -1$ for nonzero $a$,
we see that $q-3+S = 8 \,\tin{a}{b}{c}{1}{1}{1} + 4 R$.

By Hasse \cite{Hasse33}, the number of points $N$ on an elliptic curve
over $\Ff_q$ satisfies $|N - (q+1)| \le 2 \sqrt{q}$. Consider the curve
$y^2 = (x-a)(x-b)(x-c)$. The homogeneous form is $Y^2Z = (X-aZ)(X-bZ)(X-cZ)$
with a single point $(0,1,0)$ at infinity.
If $(x-a)(x-b)(x-c)$ is zero for 3 values of $x$, a nonzero square
for $m$ values of $x$, and a nonsquare for the remaining $q-3-m$ values of $x$,
then $N = 1 + 3 + 2m$ and $S = m - (q-3-m) = 2m+3-q$.
Hence $|S| = |N-(q+1)| \le 2 \sqrt{q}$.
It follows that
$\left| \tin{a}{b}{c}{1}{1}{1} - \frac{q-9}{8} \right| \le \frac14 \sqrt{q} + \frac34$.
\qed

\begin{Corollary}
If $q > 25$ then any three distinct vertices in $\Gamma$
have a common neighbor. \qed
\end{Corollary}

{\footnotesize
The table below gives for small $q$ the values of
$[h \, i \, j] := \tin{a}{b}{c}{h}{i}{j}$ that occur.
For each $q$, the first line is for triangles $abc$,
the second line for paths of length 2. The remaining cases follow by
complementation.

\smallskip
\begin{tabular}{c|cccc}
$q$ & $[1 \, 1 \, 1]$ & $[1\, 1 \, 2] $ & $[1\, 2 \, 2]$ & $[2\, 2\, 2]$ \\
\hline
5 & - & - & - & - \\
  & 0 & 0 & 2 & 0 \\
9 & 0 & 0 & 6 & 0 \\
  & 0 & 3 & 2 & 1 \\
13 & 0 & 3 & 6 & 1 \\
   & 0--1 & 3--6 & 2--5 & 1--2 \\
\end{tabular}
\quad\quad
\begin{tabular}{c|cccc}
$q$ & $[1 \, 1 \, 1]$ & $[1\, 1 \, 2] $ & $[1\, 2 \, 2]$ & $[2\, 2\, 2]$ \\
\hline
17 & 0 & 6 & 6 & 2 \\
   & 1--2 & 3--6 & 5--8 & 1--2 \\
25 & 0--2 & 6--12 & 6--12 & 2--4 \\
   & 2--3 & 6--9 & 8--11 & 2--3 \\
29 & 2 & 9 & 12 & 3 \\
   & 2--4 & 6--12 & 8--14 & 2--4 \\
\end{tabular}
\par}

\bigskip

Returning to the problem in the second paragraph, if
$\psi \colon \Ff_q^{*2} \cup \{0\} \rightarrow \Cc^*$  
satisfies conditions {\it (i)} and {\it (ii)},
then $\psi(-a)=\psi(a)^{-1}$ for each $a$ and the 
extension of $\psi$ to $\hat{\psi} \colon \Ff_q  \rightarrow \Cc^*$
via $\hat{\psi}(a+b)=\psi(a)\psi(b)$ for  $a,b \in \Ff_q^{*2}$,
is well-defined.
Given  $a,b\in \Ff_q$, we  locate $c$ with $c \adj 0,a,-b$
so that $c,a-c,b+c \in \smash{\Ff_q^{*2}}$. Now $\hat{\psi}(a+b) =
\psi(a-c) \psi(b+c) = \psi(a-c)\psi(c)\psi(-c)\psi(b+c) =
\hat{\psi}(a) \hat{\psi}(b)$,
showing for $q > 25$ that $\hat{\psi}$ is an additive character.
The cases $5 < q \le 25$ can be done by hand.

\bigskip
In the above, we gave bounds for $\tin{a}{b}{c}{1}{1}{1}$,
in particular for the number of $K_4$'s on a given triangle $abc$.
In case $q = p$ is prime, a closed formula for the total number
of $K_4$'s on a given edge was given by Evans, Pulham \& Sheehan
\cite{Evans-at-al81}. If $p = m^2+n^2$ where $n$ is odd, this number is
$\frac{1}{64}((p-9)^2-4m^2)$.

\bigskip
The bounds of Proposition \ref{bds} are best possible:

\begin{Proposition}
If $q = (4s+1)^2$ for some integer $s \ge 1$, then

(i) For a suitable triangle $abc$ one has
$\tin{a}{b}{c}{1}{1}{1}  = \frac{q-9}{8} - \frac14 \sqrt{q} - \frac34 = 2(s^2-1)$.

(ii) For a suitable cotriangle $abc$ one has
$\tin{a}{b}{c}{1}{1}{1}= \frac{q-9}{8} + \frac14 \sqrt{q} + \frac34 = 2s(s+1)$.
\end{Proposition}

\Proof
If $abc$ is a triangle or a cotriangle, then
$\tin{a}{b}{c}{1}{1}{1} + \tin{a}{b}{c}{2}{2}{2} = \frac{q-9}{4}$.
Also, $\tin{a}{b}{c}{2}{2}{2} = \tin{ea}{eb}{ec}{1}{1}{1}$
for any nonsquare $e$.
So {\it (i)} and {\it (ii)} are equivalent. Let us prove {\it (i)}, that is,
prove that $N = q - 2\sqrt{q} + 1$ occurs for a suitable curve
$y^2=(x-a)(x-b)(x-c)$ where $abc$ is a triangle.

By Waterhouse \cite{Waterhouse69} there are elliptic curves
with $N = q \pm 2\sqrt{q} + 1$ points when $q$ is a square.
A curve $y^2=(x-a)(x-b)(x-c)$ has three points of order 2, so
2-torsion subgroup $\Zz_2 \times \Zz_2$, so that its number of points
is 0 mod 4. Conversely, by Auer \& Top \cite{AuerTop02}, given an
elliptic curve $E$ with 0 mod 4 points, there is one with the same
number of points in Legendre form $y^2 = x(x-1)(x-\lambda)$,
except in case $q=r^2$ for a (possibly negative) integer
$r \equiv 1$ (mod 4) when $|E| = (r+1)^2$.
Consequently, there is a curve $y^2 = x(x-1)(x-\lambda)$ with
$N = (r-1)^2$ points. Then $S = N-(q+1) = -2r$ and
$8\,\tin{0}{1}{\lambda}{1}{1}{1} +4R = N-4 = (r-1)^2-1 = 16s^2-4$ and
$\smash{\tin{0}{1}{\lambda}{1}{1}{1}  = 2s^2 - \frac{R+1}{2}}$.
In the extreme cases, $E$ is supersingular
(e.g.~because $N \equiv 1$ (mod $p$))
and according to \cite{AuerTop02} (\S3) $\lambda$ is a square in $\Ff_{p^2}$,
and then also $1-\lambda$ is a square in $\Ff_{p^2}$,
so that $\{0,1,\lambda\}$ is a triangle and $R = 3$.
\qed

\section*{Acknowledgments}
The second author thanks Bill Kantor for helpful remarks.
The work of the second author was supported, in part,
through a grant from the National Science Foundation (DMS Award \#1808376)
which is gratefully acknowledged.

\end{document}